\documentclass[preprint,1p,12pt]{elsarticle}

\usepackage{hyperref}

\usepackage{amsmath, amssymb, mathrsfs}
\usepackage{amsthm}
\usepackage{graphicx,graphics,color}
\usepackage{lineno}
\usepackage{natbib}
\usepackage{ulem}
\usepackage{amsmath}
\usepackage{amssymb}
\usepackage{fancyhdr}

\begin{document}

\begin{frontmatter}

\title{A semi-decentralized control strategy for urban traffic}

\author[add1]{Nadir Farhi\corref{cor1}}
\ead{nadir.farhi@ifsttar.fr}
\cortext[cor1]{Corresponding author}

\author[add1]{Cyril Nguyen Van Phu}

\author[add1,add2]{Mouna Amir}

\author[add1]{Habib Haj-Salem}

\author[add1]{Jean-Patrick Lebacque}

\address[add1]{Universit\'e Paris-Est, IFSTTAR / COSYS / GRETTIA, F-77447 Marne-la-Vall\'ee cedex 2, France.}
\address[add2]{Universit\'e de Versailles Saint-Quentin-en-Yvelines, France.}

\begin{abstract}
We present in this article a semi-decentralized approach for urban traffic control, based on the
TUC (Traffic responsive Urban Control) strategy. We assume that the control is centralized as in
the TUC strategy, but we introduce a contention time window inside the cycle time, where
antagonistic stages alternate a priority rule. The priority rule is set by applying green colours
for given stages and yellow colours for antagonistic ones, in such a way that the stages with green
colour have priority over the ones with yellow colour. The idea of introducing this time window is
to reduce the red time inside the cycle, and by that, increase the capacity of the network junctions.
In practice, the priority rule could be applied using vehicle-to-vehicle (v2v) or
vehicle-to-infrastructure (v2i) communications. The vehicles having the priority pass almost normally
through the junction, while the others reduce their speed and yield the way. We propose a model for
the dynamics and the control of such a system. The model is still formulated as a linear quadratic
problem, for which the feedback control law is calculated off-line, and applied in real time. The
model is implemented using the Simulation of Urban MObility (SUMO) tool in a small regular
(American-like) network configuration. The results are presented and compared to the classical TUC strategy.
\end{abstract}


\end{frontmatter}

\section{Introduction}
\label{introduction}


Recent advances in information and communication technologies improve vehicular traffic in urban road
networks by enabling the development of innovative urban traffic control strategies.
While the traffic control in urban road networks is still done by setting traffic lights, intelligent
transportation systems (ITS) are being tested in many cities.
Various agents in the road network will be able to communicate from vehicle to vehicle (V2V) or
from vehicle to infrastructure (V2I) for example. Big data sets, with different levels of information
(microscopic, macroscopic) will be processed in real time and adaptive control strategies will be applied.
The whole process of urban traffic control needs to be redefined in order to take into account this development. 

Several levels of information are distinguished in the big amounts of data that are made available by ITS.
The whole information cannot be optimally exploited with a unique centralized or distributed traffic control system.
A multi-level control system needs to be developed in order to optimally use each level of information for the corresponding
control level. Macroscopic information could be transmitted to the centralized controller, while the microscopic one could be used by
the local controller, which should operate in a short time horizon, compared to the high-level controller.
Multi-level control schemes have been recently proposed; see for example~\citep{RHG15, Var13}. In~\citep{RHG15}, the
control uses macroscopic fundamental diagrams (MFD)~\citep{GD07,DG08b,FGQ05,FGQ07,Far08,Far09,FGQ11}.


Using traffic lights, the main urban traffic parameters are: phase specification, split, cycle time, and offset. 
Fixed time urban traffic control (UTC) strategies appeared in the 1950s with coordination of signals that optimizes the offsets.
These strategies use historic datasets, and therefore, are unable to adjust to changing conditions.
The most well-developed and widely used UTC system is TRANSYT~\citep{Rob69}.
With advances in detection, communication, data processing, and control strategies, traffic responsive UTC systems appeared, where
centralized and distributed systems are distinguished.
Among the main centralized ones, we cite SCOOT \citep{HRBW81,BWB98}, SCATS~\citep{Low82}, RHODES~\citep{HMS92}, MOTION~\citep{Bus96}, and 
TUC~\citep{Dia99}.
For distributed responsive UTC, we cite UTOPIA~\citep{DMRV84}, PRODYN~\citep{FKL90}, OPAC~\citep{Gar91}.
Other UTC systems define an intermediate level of centralization.


Traffic responsive UTC systems use feedback controls on the state of the traffic and permit, by that, to meet traffic demand.
Moreover, the control may be set in such a way to be robust, in the sense that it responds rapidly to
disruptions. Furthermore, such controls are automatically adaptive to works and operations, and so installation and maintenance costs
are reduced.
    

We propose in this article an extension of the traffic responsive urban control strategy 
TUC (Traffic Urban Control)~\citep{Dia99, DPA02, DDAPBSL03}.
Our extension introduces a kind of decentralization in the optimization of the right of way assignment.
We introduce a contention time window inside the cycle time, where
antagonistic stages alternate a priority rule. The priority rule is set by applying green colours
for given stages and yellow colours for antagonistic ones, in such a way that the stages with green
colour have priority over the ones with yellow colour.
A TUC-based centralized control determines the optimal split of green, 
red and yellow lights at the level of every junction. A decentralized
system manages the traffic of antagonistic stages during the yellow signal, taking into account the characteristics
of each junction.
By doing this, we aim to reduce the red time inside the cycle, increase the capacity of the network, and 
reduce users' delays.
The traffic management during the yellow times would be realized based on vehicle to vehicle (v2v) and/or vehicle to
infrastructure (v2i) communications.

We present in this article preliminary results of this semi-decentralization on a small American-like city.
The results demonstrate the efficiency of this extension with respect to the classical TUC control.
On a selected scenario of traffic demand, we show that the semi-decentralized TUC controls better the traffic, in
the sense that it is able to respond efficiently and rapidly to congestion.

\section{A short review of TUC}
\label{review}

\begin{figure}[htbp]
  \begin{center}
    \includegraphics[width=6cm]{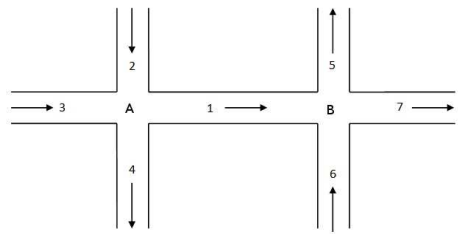}    
    \caption{Academic example explaining the TUC strategy.}
    \label{fig1}
  \end{center}
\end{figure}

TUC~\citep{Dia99, DPA02, DDAPBSL03} is a coordinated control strategy based on a store-and-forward approach.
It can be implemented for large-scale networks, in real time, even under saturated traffic conditions.
The split control part of TUC varies the green-stage durations of all stages at all the junctions of a urban network
around given nominal values, and under a simplified traffic dynamics. 
The objective is to avoid oversaturations and spillbacks of link queues.
In order to briefly explain the approach, let us consider the small network of Figure~\ref{fig1}, with the following 
notations.

\begin{tabular}{ll}
$c$			& cycle time duration, in seconds.\\
$k$			& discrete time index, corresponding to a duration of $kc$ sec.\\
$x_i(k)$		& number of cars on link $i$ at discrete time $k$.\\
$\bar{x}_i$		& constant nominal number of cars on link $i$.\\
$\Delta x_i(k)$		& $= x_i(k) - \bar{x}_i$.\\
$s_i$			& saturation flow on link $i$.\\
$g_i(k)$		& green time duration for link $i$ during the $k$th cycle.\\
$\bar{g}_i$ 		& constant nominal green time duration for the stream coming \\
                        & from link $i$. \\
$\Delta g_i(k)$ 	& $= g_i(k) - \bar{g}_i$. \\
$u_i(k)$		& $= (g_i(k)/c) s_i$ average outflow from link $i$ during the $k$th cycle.\\
$d_i(k)$ 		& arrival demand flow to link $i$ at discrete time $k$. \\
$\bar{d}_i$ 		& constant nominal arrival demand flow to link $i$. \\
$\Delta d_i(k)$ 	& $= d_i(k) - \bar{d}_i$. \\
$\alpha_{ij}$		& turning movement ratio from link $i$ to link $j$ \\
\end{tabular}
\newline
\newline

The definition of $u_i(k)$ assumes sufficient demand on link $i$. Note that the oscillations of vehicle queues in the links due to green/red
communications, and the effect of offset for consecutive junctions cannot be described by the model.

According to Figure~\ref{fig1}, the number of cars on link 1 is updated as follows.

\begin{equation} \label{dyn0}
  x_1(k+1) = x_1(k) + d_1(k) + \alpha_{21} s_2 g_2(k) + \alpha_{31} s_3 g_3(k) - s_1 g_1(k).
\end{equation}

Then, by introducing the nominal amounts, and by using vectorial notations, we get :

\begin{equation}\label{dyn1}
  \Delta x(k+1) = \Delta x(k) + B \Delta g(k) + D \Delta d(k),
\end{equation}
where $B$ and $D$ are matrices built basing on the dynamics~(\ref{dyn0}) written on the whole network.

Assuming that the variations of the arrival demand flows on every link inside the cycle time sum to zero, we get the following 
linear system :

\begin{equation}\label{dyn2}
  \Delta x(k+1) = \Delta x(k) + B \Delta g(k),
\end{equation}

Bounds for minimum green times and maximum storage capacity of links must also be considered.

The criterion is the following, where $\lambda$ is a discount factor, and where an infinite time horizon is considered.
\begin{equation}\label{crit1}
  J = \min_{\Delta g} \frac{1}{2} \sum_{k=0}^{+\infty} \frac{1}{(1+\lambda)^k} \left( \|\Delta x(k)\|_Q^2 + \|\Delta g(k) \|_R^2\right),
\end{equation}
where $Q$ and $R$ are non-negative definite, diagonal weighting matrices.
The first term on~(\ref{crit1}) aims to minimize the risk of oversaturation and the spillback of link queues, while the second term
is used to influence the magnitude of the control.

The control bounds are treated externally of the LQ problem solving.
The solution for such problems consists in solving an algebraic Riccati equation, which then leads
to the following optimum feedback control :
\begin{equation}\label{feed1}
  g(k) = \bar{g} - L x(k).
\end{equation}
where $L$ is the gain matrix; see~\citep{Dia99, DPA02, DDAPBSL03} for more details.

\section{Semi-decentralization}

The model we present here is an extension of the classical model presented above. Instead of considering only green and red time durations
in a cycle time (in addition to the lost time, which we consider implicit here and for which we assign the orange colour), 
we also consider yellow time durations.
The objective here is to reduce the red time duration. To do that, we divide this duration into two time periods : red and yellow.
By that, when a stage is assigned a red or a yellow time, the antagonistic stage is assigned a green time.

\begin{figure}[htbp]
  \begin{center}
    \includegraphics[width=10cm]{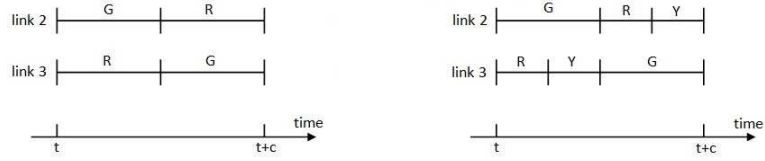}    
    \caption{The cycle time in the classical model, and in the new model. G: green, R: red, Y: yellow.}
    \label{cycle-time}
  \end{center}
\end{figure}

We notice here that our model is an extension of the classical TUC model, because it is sufficient to set the yellow times to zero to get the 
classical model.

In order to explain the model, let us consider junction $A$ of example of Figure~\ref{fig1}. 
Only two stages can be considered here, each of them with only one stream. 
One stage is associated to link~2  and the other to link~3.
In this case, and in the classical TUC model, at every cycle $k$, we only have one independent control variable on that junction,
which is the green or red duration
of any of the two streams. All the other time durations are dependent variables.
We consider $g_2(k)$: the green time duration for link 2 as the independent control variable, then the dependent variables
can be easily obtained as follows :
\begin{itemize}
 \item $r_2(k) = c - g_2(k)$ : red duration for link 2
 \item $g_3(k) = r_2(k)$ : green duration for link 3
 \item $r_3(k) = g_2(k)$ : red duration for link 3
\end{itemize}

By considering yellow time durations, we  need to choose three independent control variables, among six variables.
For example the following three independent control variables can be considered.
\begin{itemize}
 \item $g_2(k)$ : green time duration for link 2
 \item $y_2(k)$ : yellow time duration for link 2
 \item $y_3(k)$ : yellow time duration for link 3
\end{itemize}
The other three dependent control variables are given as follows (see figure~2): 
\begin{itemize}
     \item $r_2(k) = c - g_2(k) - y_2(k)$ : red time duration for link 2
     \item $r_3(k) = g_2(k) - y_3(k)$ : red time duration for link 3
     \item $g_3(k) = c - g_2(k)$ : green time duration for link 3
\end{itemize}

\subsection{The dynamics}

Let us consider the following additional notations.
\begin{itemize}
  \item $q^{\max}_{J}$ : capacity (maximum flow) of junction $J$.
  \item $Q_{ij}(k)$ : total flow going from link $i$ to link $j$ during the $k$th cycle.
  \item $Q^{out}_i$ : total flow exiting from link $i$ during the $k$th cycle.
  \item $\gamma_J$ : friction coefficient on junction $J$, with $0\leq \gamma \leq 1$.
\end{itemize}

We write the traffic dynamics on link 1 of Figure~\ref{fig1} with the new control model.

\begin{equation} \nonumber
  x_1(k+1) = x_1(k) + d_1(k) + Q_{21}(k) + Q_{31}(k) - Q^{out}_1(k),
\end{equation}
with $Q_{21}(k), Q_{31}(k)$ and $Q^{out}_1(k)$ are given in~(\ref{eqq1})-(\ref{eqq3}), where we introduce a new parameter $\gamma_J$ 
(for junction $J$) which we call here a friction coefficient, and which expresses the 
bother between vehicles entering into the junction from antagonistic stages during the contention time window.
For example, in~(\ref{eqq1}), the flows of vehicles going from link 2 to link 1 during different time durations of 
the $k$th cycle are given as follows.
\begin{itemize}
 \item During $r_3(k) = g_2(k) - y_3(k)$, the flow is $\alpha_{21}s_2(g_2(k)-y_3(k))$, as usual.
 \item During $y_3(k)$, the flow is $\alpha_{21}s_2y_3(k)$ as usual, but multiplied by the friction coefficient
   $\gamma_A$ between the streams coming from link 2 (with green time) and link 3 (with yellow time),
   since the local control is activated with a priority rule setting. Link 2 has priority over link 3 during this time period.
 \item During $r_2(k)$, the flow is zero.
 \item During $y_2(k)$, the stream coming from link 3 has priority over the one coming from link 2.
   Therefore, the whole junction capacity $q_A^{\max}y_2(k)$ is used by the stream of link 3, and
   the remaining capacity $q_A^{\max}y_2(k) - s_3 y_2(k)$ is used by link 2. This flow is also multiplied
   by the coefficient friction $\gamma_A$ since the two streams pass through junction $A$ during the same time period.
\end{itemize}

\begin{align}
   Q_{21}(k) & = \alpha_{21} s_2 (g_2(k) - y_3(k)) \nonumber \\
             & + \gamma_A \alpha_{21} s_2 y_3(k) + \gamma_A (q_A^{\max} y_2(k) - s_3 y_2(k)). \label{eqq1} \\
   Q_{31}(k) & =  \alpha_{31} s_3 (c - g_2(k) - y_2(k)) \nonumber \\
             & + \gamma_A \alpha_{31} s_3 y_2(k) + \gamma_A (q_A^{\max} y_3(k) - s_2 y_3(k)). \label{eqq2} \\
   Q^{out}_1(k) & = s_1 (g_1(k) - y_6(k)) \nonumber \\
                & + \gamma_B s_1 y_6(k) + \gamma_B (q_B^{\max} y_1(k) - s_6 y_1(k)). \label{eqq3}
\end{align}

The dynamics~(\ref{eqq1})-(\ref{eqq3}) are still linear on the variables $x_i, g_i$ and $y_i$.
We notice here that the dynamics are written with only independent controls. As it has already been explained above,
on junction $A$, for example, the independent controls are $g_2, y_2$ and $y_3$.
As in the classical TUC model, we consider nominal demands $\bar{d}_i$, nominal numbers of cars $\bar{x}_i$ and nominal
independent controls $\bar{g}_i$ and $\bar{y}_i$.
The choices of $\bar{x}$ and $\bar{g}$ can be done by the same way as in the classical TUC model. 
One way to choose $\bar{y}$ is to take  $\bar{y}_i = c - \bar{g}_i$. This is equivalent to say that the nominal red time is zero.
This choice can also be dependent on the junction design.
Then it is very easy to derive a linear dynamics similar to~(\ref{dyn2}).
For the criterion we take exactly the one of~(\ref{crit1}), written with the new (independent) control variables $\Delta g_i$.
Again, a linear quadratic problem is obtained, and the optimal control is derived by solving a Riccati equation as in the classical TUC model. 

\section{Numerical example}

\begin{figure}[htbp]
  \begin{center}
    \includegraphics[width=10cm]{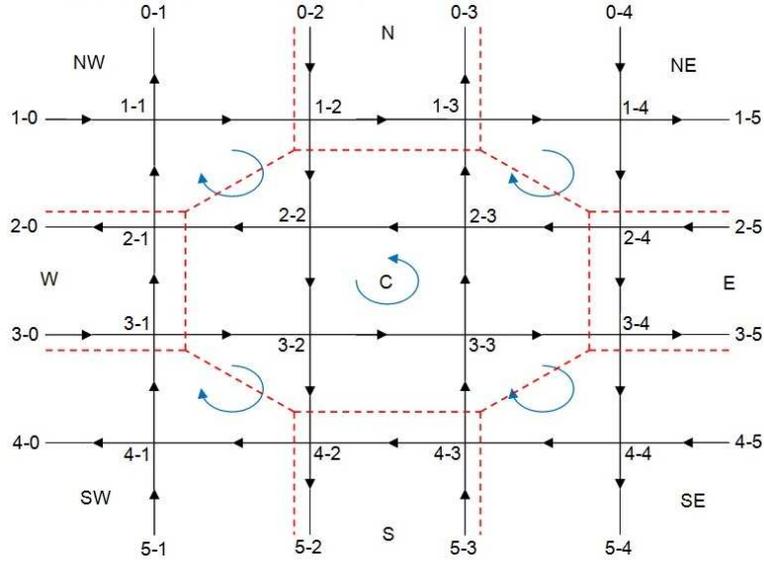}    
    \caption{Regular network example.}
    \label{example1}
  \end{center}
\end{figure}

In this section, we apply the control model presented above, on a small regular (American-like) network
of four horizontal and four vertical roads, with alternated directions, as shown in Figure~\ref{example1}.

For the saturation flow values, we take the recommended ones in urban networks ($s_i = 1800 veh./h, \forall i$, as shown in Table~\ref{tab2}),
without corrective factors; see for example~\citep{Coh93}.
To compute the optimal cycle, we consider here a fixed cycle time that we approximate to 60 seconds, using the Webster Method~\citep{Web58}: 
$c = (1.5 T + 5)/(1 - Y)$, where $T$ is the total lost time per cycle, $Y$ is the junction load.
The cycle time is then projected onto the interval $[40 s, 90 s]$.

\subsection{Model implementation and Simulation Tools}

We used SUMO, see for example \citep{behrisch2011sumo}, and its interface TRACI \citep{TraCI} to simulate and implement the model.
The source code has been written in Python. The main tasks were : 
\begin{itemize}
 \item  build the network topology and the demand using SUMO tools and original configuration files.
 \item  design an algorithm and the source code architecture that enable the construction of the $B$ matrix in equation~(\ref{dyn2}).
 \item  implement the contention time window and the associated priority rule.
 \item  solve a Riccati equation, and at every cycle, measure the state, and apply the control on the traffic light signals.
 \item  analyze the simulation data outputs, including state and control vectors, by rendering graphical results.
\end{itemize}

The time contention window is implemented as follows.
On a given junction, and inside such contention time window, we consider first vehicles in incoming approaches.
We compute the distances from those vehicles to the junction.
In order to avoid conflicts, 
at a given time in the time window, if the distance to the junction, of the first vehicle on the link with yellow stage, 
is less than a given constant distance $m$, and if the distance to the junction, of the first vehicle on a link with green antagonistic stage, 
is less than a given constant distance $M$, we slow down the first vehicle on the link with yellow stage.

In general, the vehicles moving on an approach with a green stage (priority approach) pass through the junctions without
checking for the antagonistic approaches. However, the vehicles moving on the approaches with yellow stages slow down at 
a distance $m$ to the junction, to check if there is any vehicle coming from an antagonistic approach with green stage.

In the numerical example we consider in this article, we chose $m=15$ meters and $M=50$ meters.
Our choice takes into account the reaction time of the drivers in SUMO, and also the simulation step length.

We plan to implement this conflict management using a communication simulator, for example the Network Simulator~\citep{mccanne1997network}.

\subsection{Network configuration}

We discuss here, the configuration of the network of Figure 3.
In this network, circuits are formed.
We distinguish two types of circuits. The central circuit in which vehicles turn 
in the anticlockwise direction, and the other four circuits in which vehicles move in the clockwise direction. 
As already shown in \citep{Far08, FGQ11}, the car-densities on the circuits of links are determinant in the stage transition
of a vehicular network. Indeed, if a circuit is full of vehicles, then a deadlock occurs and spreads on the network.

In the network we consider here, the central circuit (which we call here the main circuit) is critical compared to the other four
circuits, (which we call here the secondary circuits).
Indeed, the secondary circuits have exits that are not constrained by any output supply, and they are closer to the borders.

In case of congestion, we need to clear out vehicles from the main circuit in order to improve the traffic,
so that the number of vehicles we take out is bigger than the one we take into the circuit.
Hence, for that circuit, the controller needs to favour the vehicles coming from the left side at the level of the four junctions
around the main circuit.
For example, if we take symmetric turn ratios, half of vehicles leaving the approaches are likely to leave the circuits,
while the other half of vehicles are likely to remain in the circuit. 
However, when the way is given to the vehicles coming from the right (with respect to the junction), half of those vehicles are likely
to enter to the circuit, while the other half is likely to not enter to the circuit.
For the secondary circuits, in case of congestion, the control shall favour vehicles coming from the right side links at the level
of the junctions associated to those circuits, in order to clear them out.

The four junctions of the main circuit are shared with other secondary circuits.
We think that the control needs to foster the evacuation of the main circuit with respect to the secondary circuits.
Therefore, the control should favour the vehicles coming from the left side approaches to the main circuit.

\subsection{Preliminary results}

We present in this section the preliminary results we obtained.
For the traffic demand, we took the scenario of Table~1.
In this scenario, we have some traffic demand inside the network. This permits us to attain saturated and congested stages.
In the other side, the traffic demand from and towards the central zone is low comparing to that from and towards the boundary zones.
This choice makes the states of the traffic controllable in the central zone of the network.

\begin{table}[h]
\caption{The traffic demand.}
\begin{tabular*}{\hsize}{@{\extracolsep{\fill}}lll@{}}
\hline
 & Central zone & Other zones\\
\hline
Central zone &   0 &  40 ({\it{veh / h}})\\
Other zones  & 40 ({\it{veh / h}}) &  250 ({\it{veh / h}})\\
\hline
\end{tabular*}
\end{table}

The other parameters are given in Table 2, where
\begin{itemize}
  \item $r$ is a positive scalar such that $Q = I$ and $R=r I$, with $I$ the associated identity matrix,
  \item $g_{i-min}$ is the minimal green time duration on link $i$,
  \item $l_i$ is the length of link $i$.
\end{itemize}  

\begin{table}[h]
\caption{The values of other parameters.}
\begin{tabular}{ccccccccc}
\hline
 $r$ & $\lambda$ & $\bar{x}_i$ & $s_i$ & $\bar{g}_i$ & $g_{i-min}$ & $c$ & $l_i$ & $\alpha_{ij}$ \\
\hline
 $0.5$ & $0.1$ & $10.5$ veh/m & $1800$ veh/h & $30$ s & $4$ s & $60$ s & $300$ m & 0.5 \\ 
\hline
\end{tabular}
\label{tab2}
\end{table}

\begin{figure} \label{fig-sim}
  \begin{center}
    \begin{tabular}{c||c}
      \includegraphics[width=5cm]{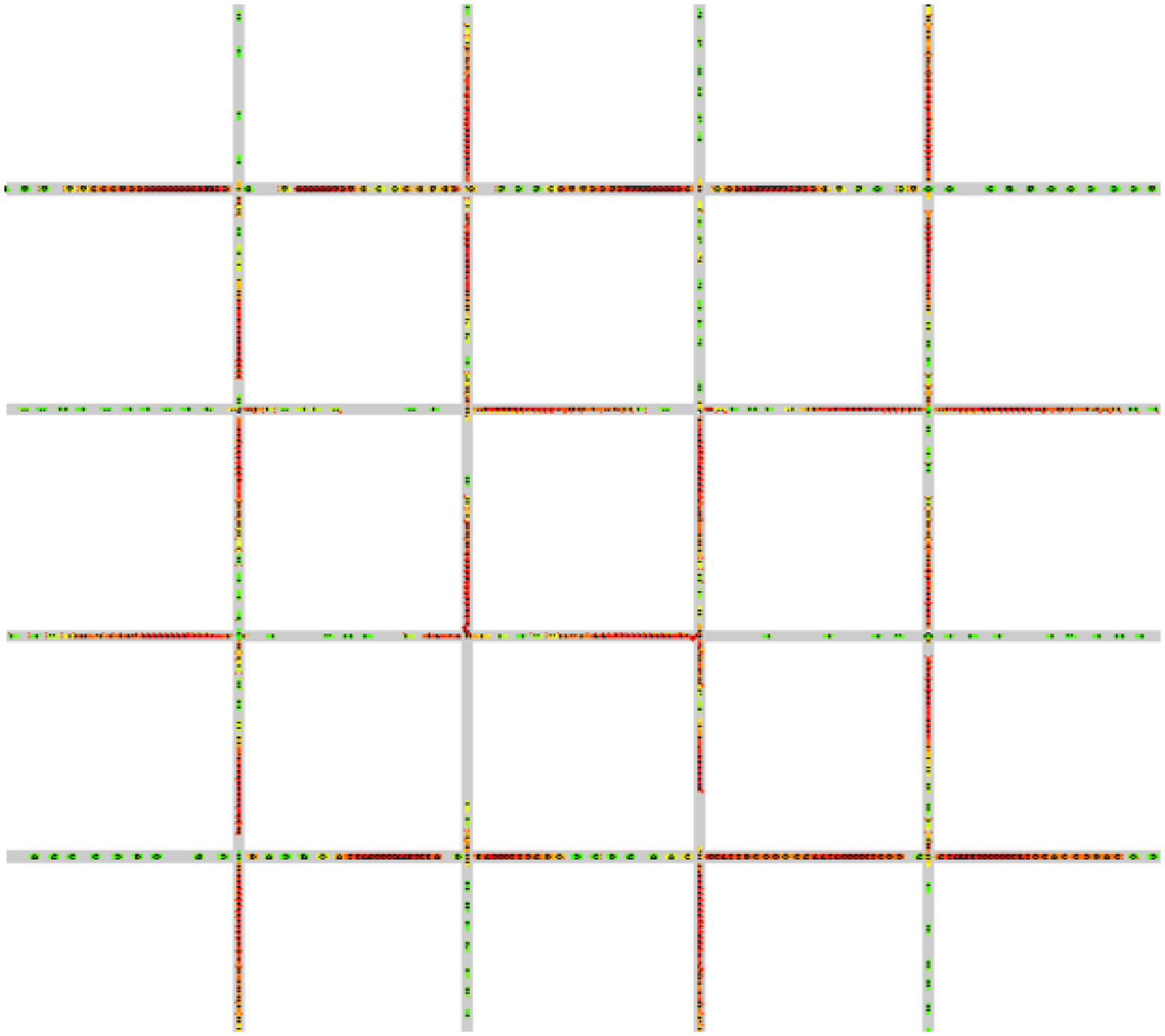} \hspace{1cm} & \hspace{1cm} \includegraphics[width=5cm]{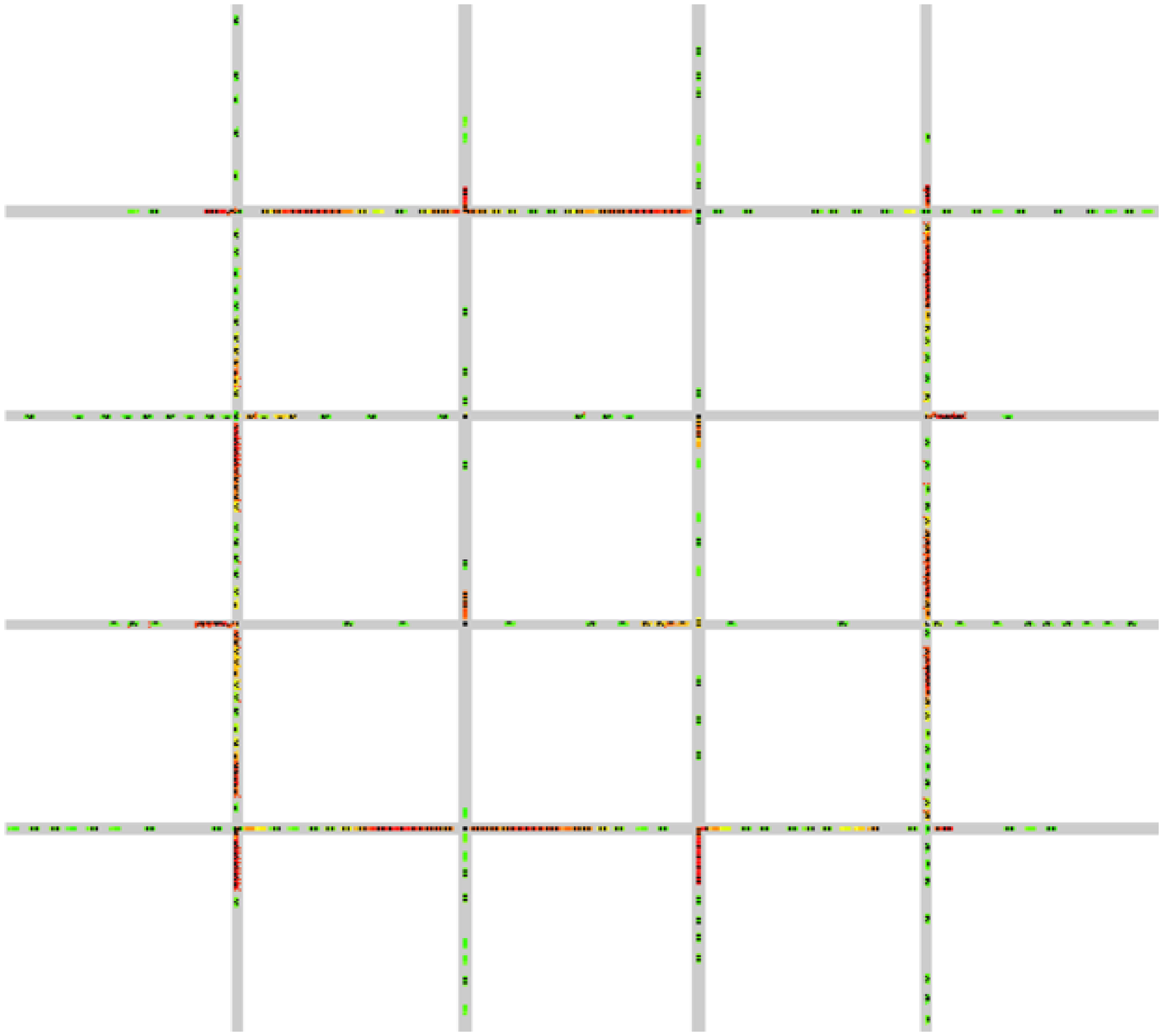}
    \end{tabular}
    \caption{The state of the traffic at the end of simulation. The colours of vehicles correspond to their speed
      (green: high speed, red: low speed).
      On the left side: Centralized TUC. On the right side: semi-centralized TUC.}
    \label{figNet}
  \end{center}
\end{figure}

\begin{figure}[htbp]
\begin{center}
\begin{tabular}{c}
  \includegraphics[width=8cm]{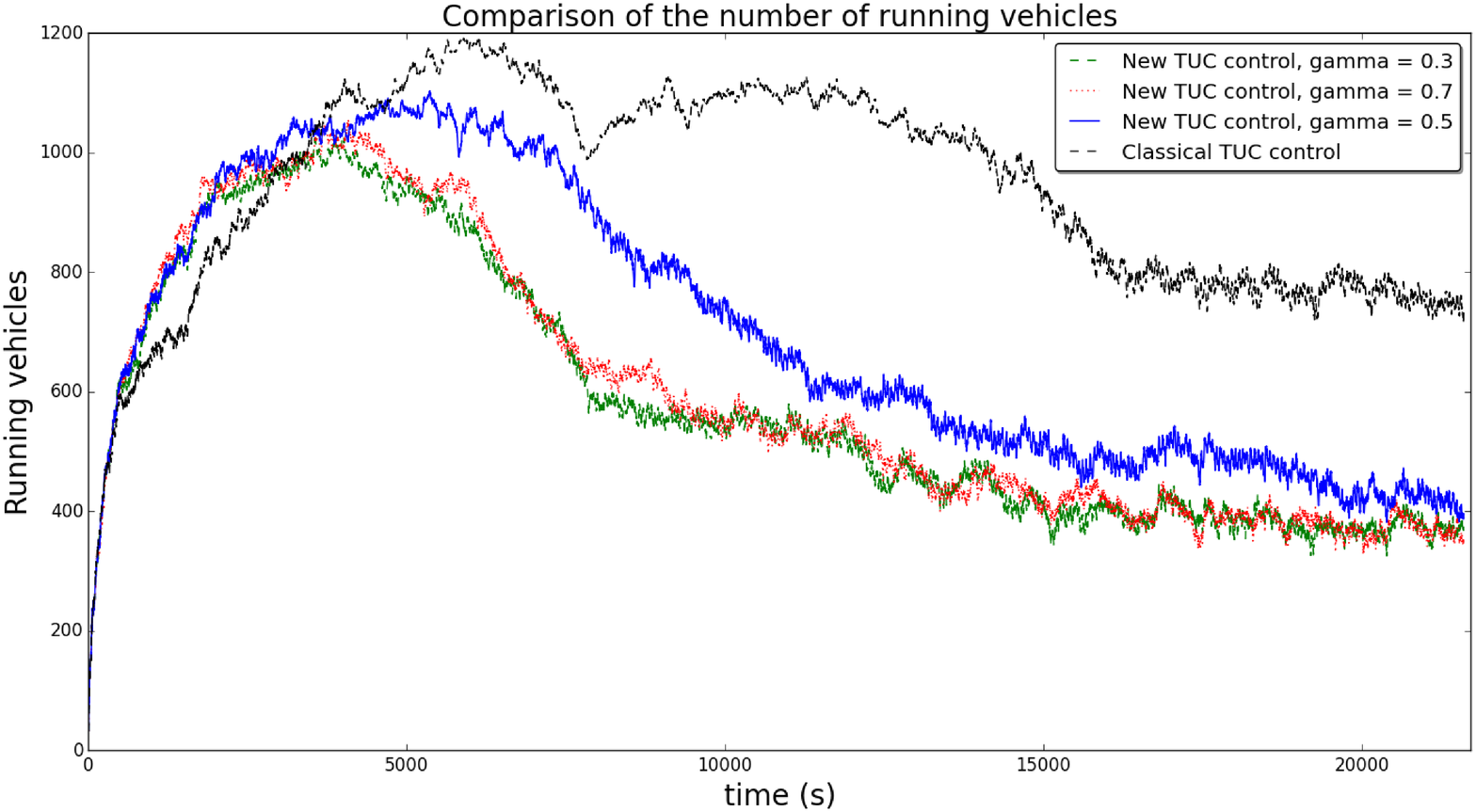} \\
  \includegraphics[width=8cm]{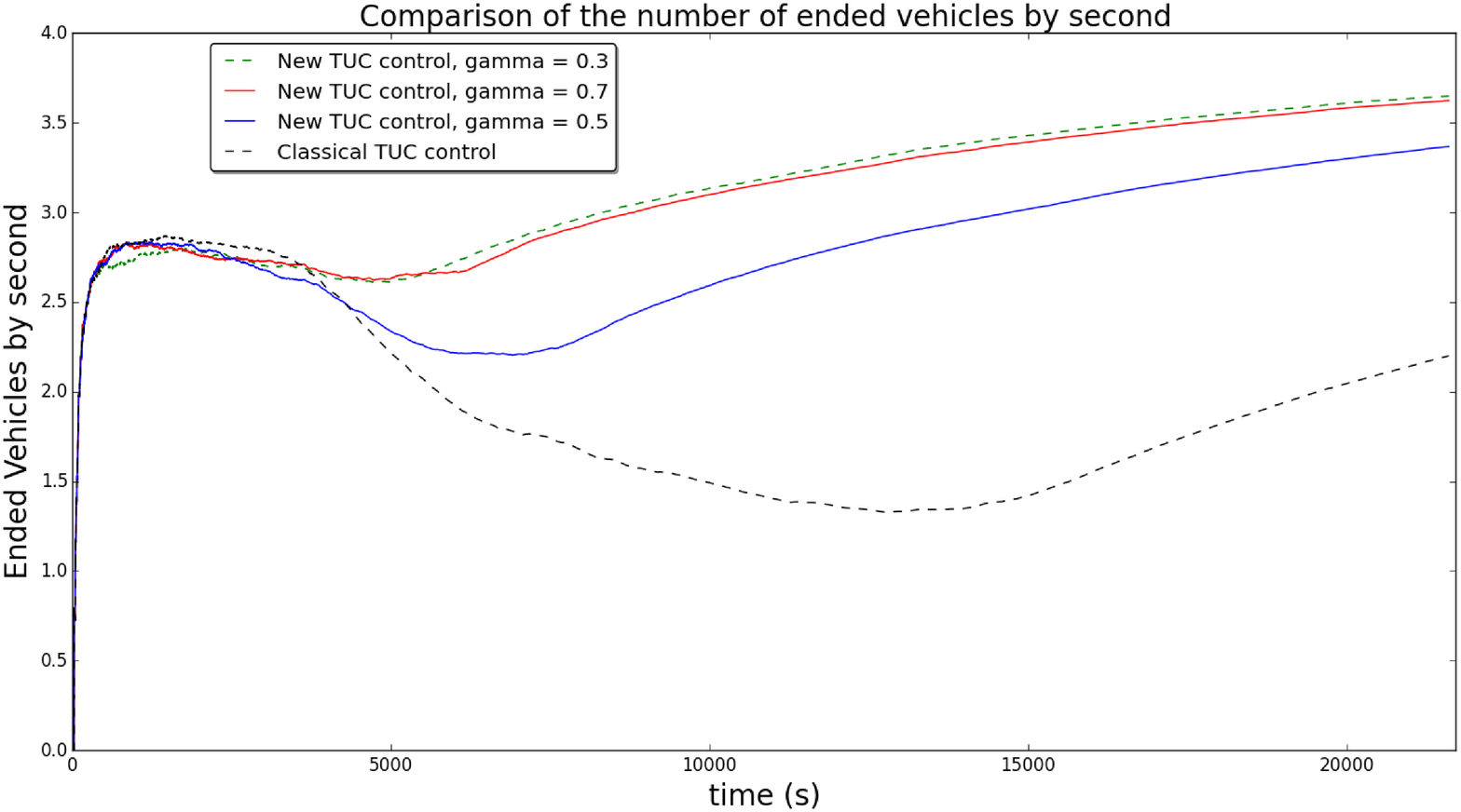}\\
  \includegraphics[width=8cm]{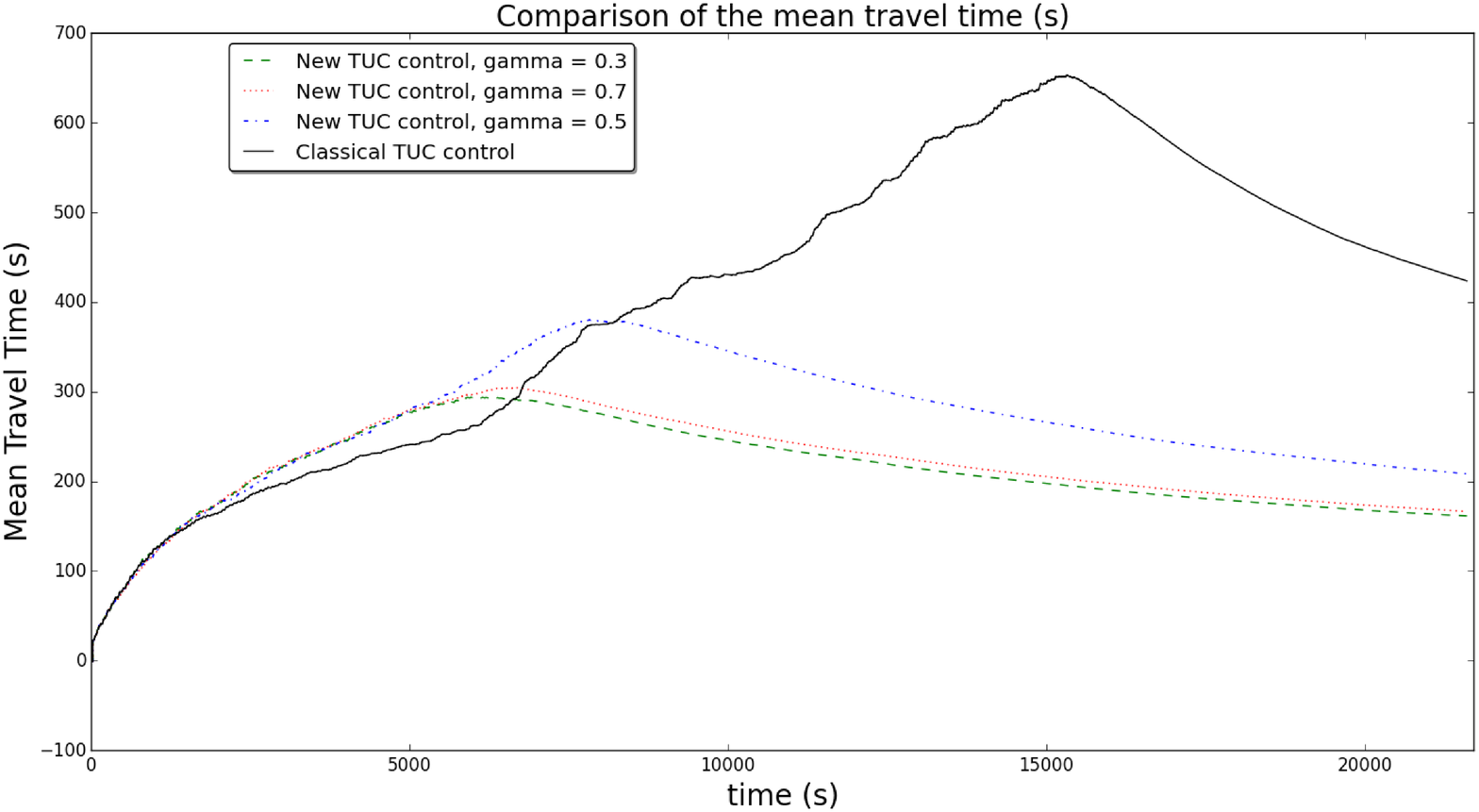}
\end{tabular}
\caption{Comparison of the classical TUC with the semi-decentralized TUC in terms of the number of running vehicles 
on the network, the flow of ended vehicles, and the average travel time through the network, respectively, in function
of time.}
\label{fig-running}
\end{center}
\end{figure}

In Figure~\ref{figNet}, we give the state of the traffic at the final time of simulation.
The evolution over time of the running vehicles in the network is given on Figure~\ref{fig-running}, where
we compare the classical TUC control with our semi-decentralized control by varying the value of the friction
parameter $\gamma$ in $\{0.3, 0.5, 0.7\}$.
We see that with our semi-decentralized control, the car-density is limited in order to optimize the capacity
of the network. The best result is obtained with $\gamma = 0.3$.

In Figure~\ref{fig-running}, we also compare the two controls in term of the cumulated ended cars through the time, and in term of
the average travel time of cars in the network. We see clearly that our control improves the whole capacity of the network.
Indeed a congestion appeared at a time around 1000 seconds. We observe that as long as the simulation runs, the two controls clear the
congestion, but the semi-decentralized control do it very rapidly compared to the centralized one. We see clearly that the difference
between the number of running vehicles decreases over time, but, even at the final time of simulation (which is 6 hours here), this difference
is still important. Figure~\ref{figNet} tells clearly that the state of the traffic with the two controls is different (fluid with
the decentralized control, and saturated with the centralized one). These results are confirmed by Figure~\ref{fig-running}, where we compare
the running and the ended vehicles, as well as the average travel time of the cars through the network.

\begin{figure}[htbp]
\begin{center}
\begin{tabular}{cc}
  \includegraphics[width=65mm]{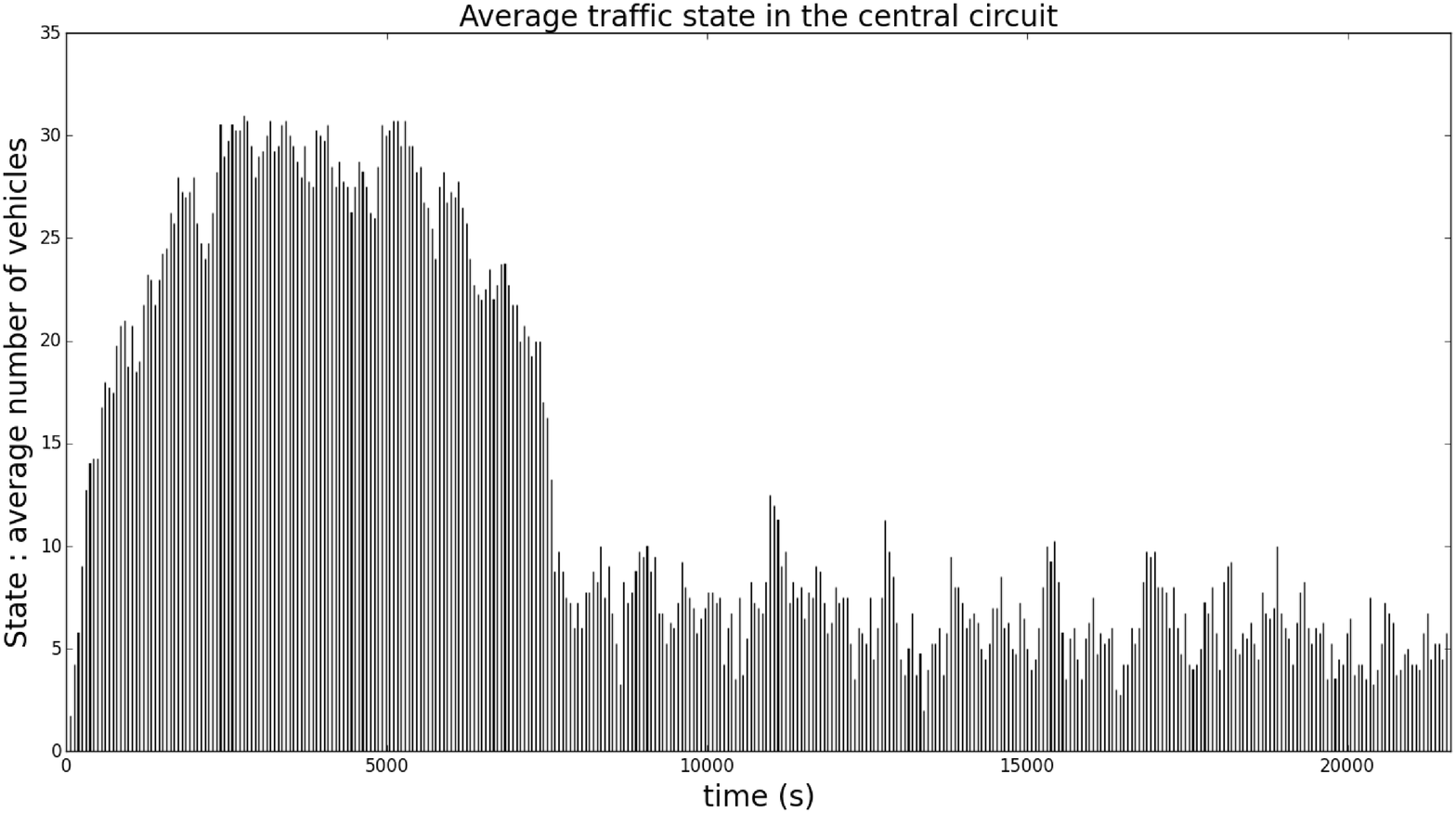} & \includegraphics[width=65mm]{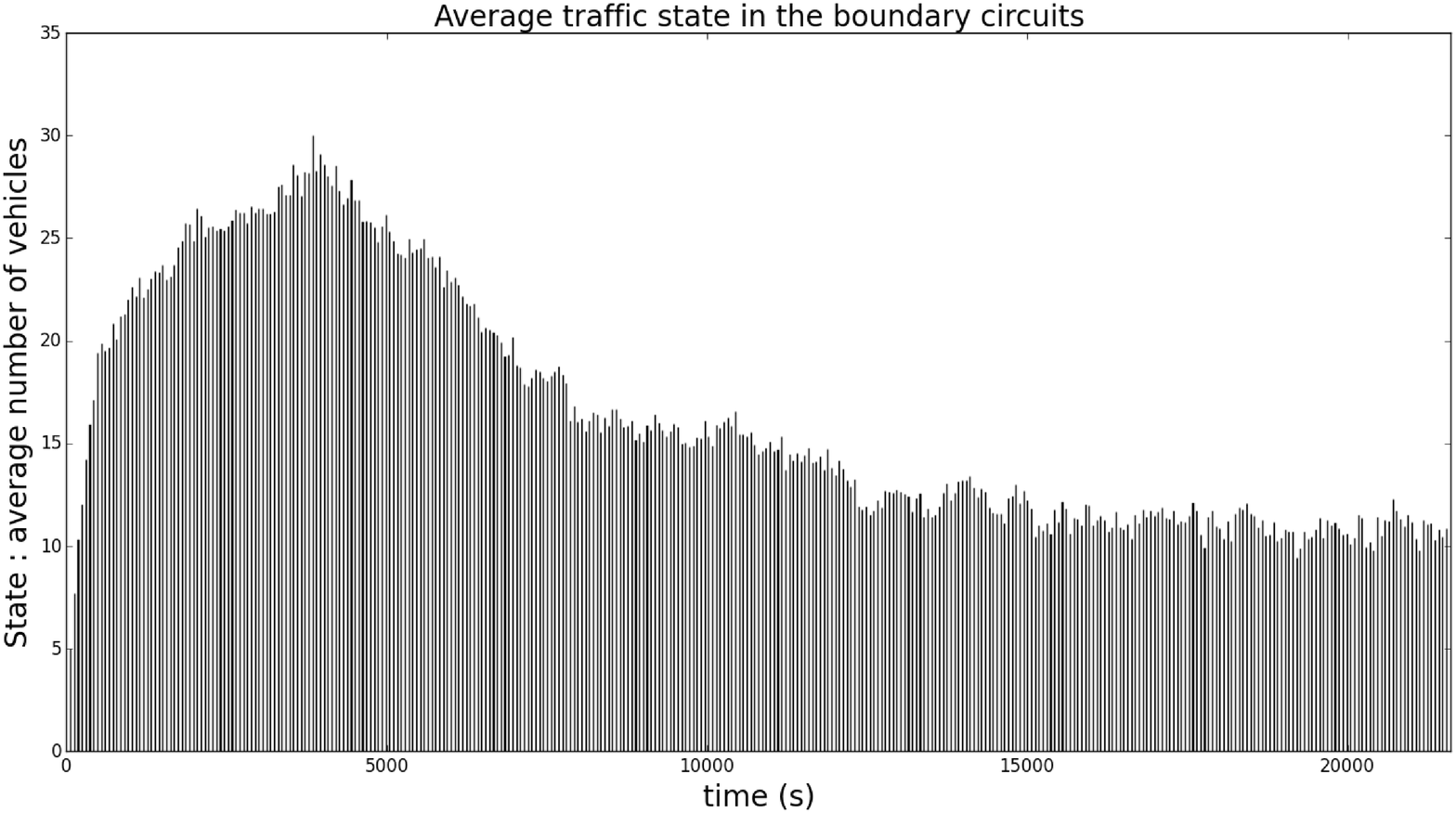} \\
  \includegraphics[width=65mm]{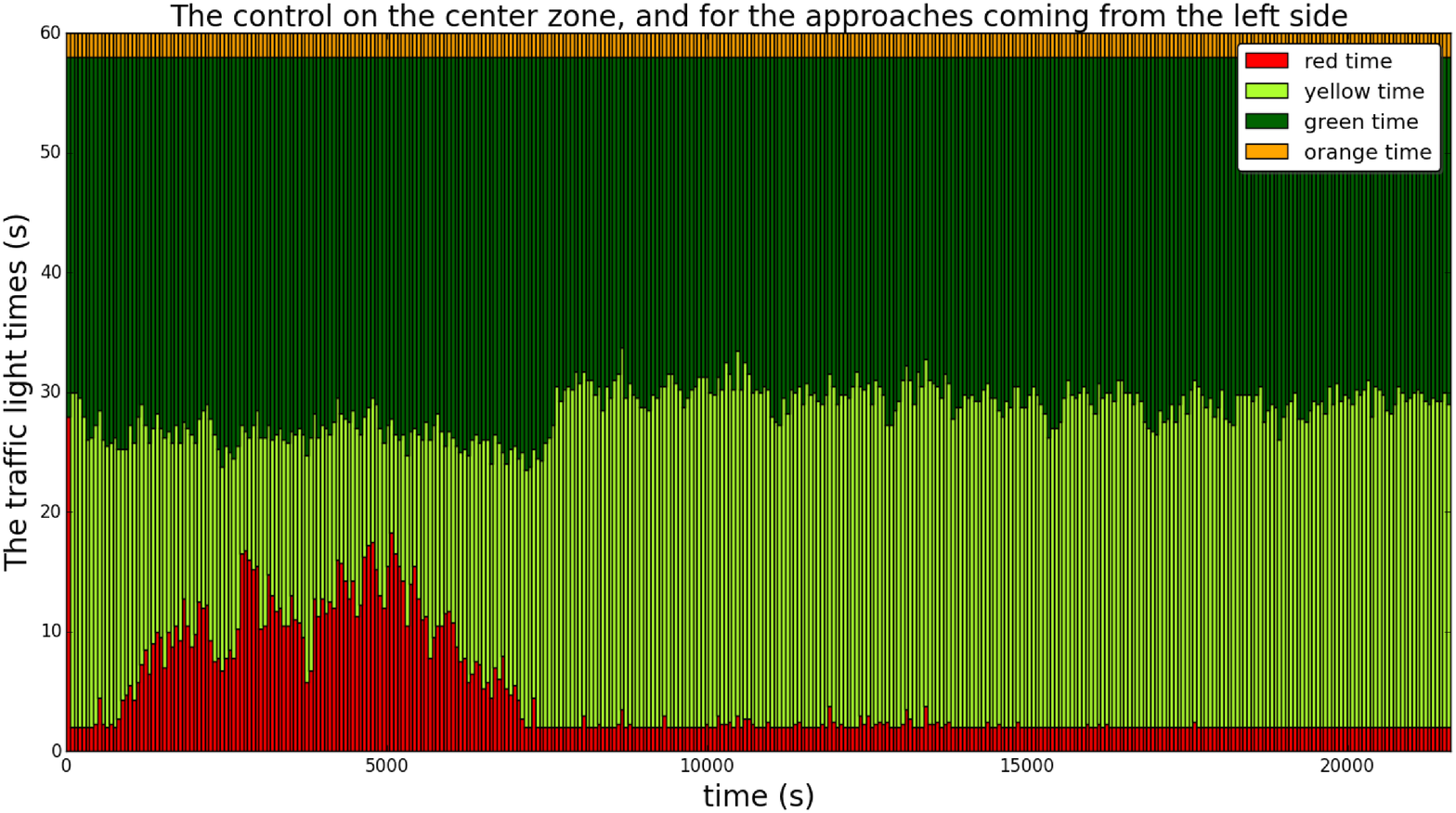} & \includegraphics[width=65mm]{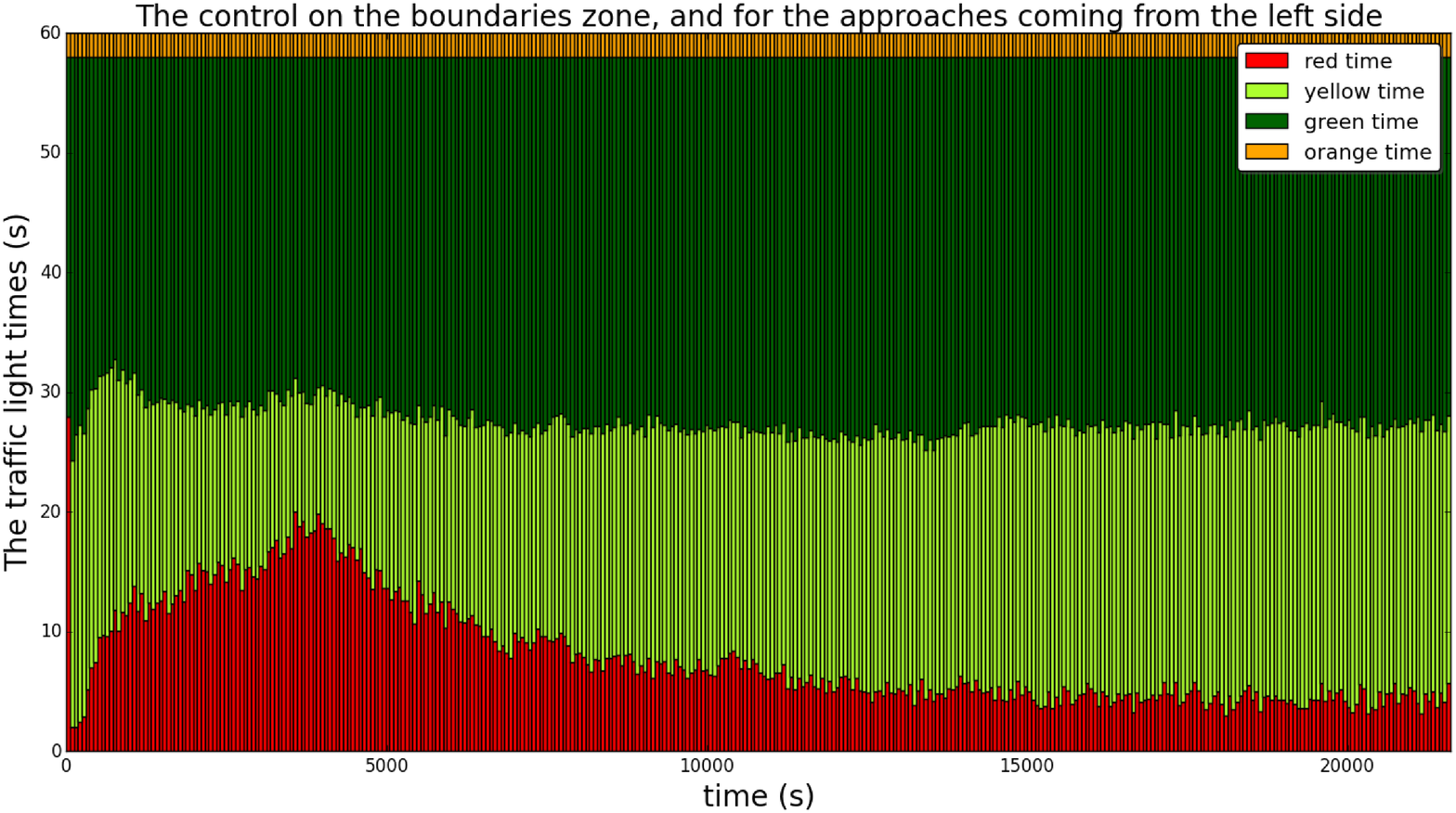} \\
  \includegraphics[width=65mm]{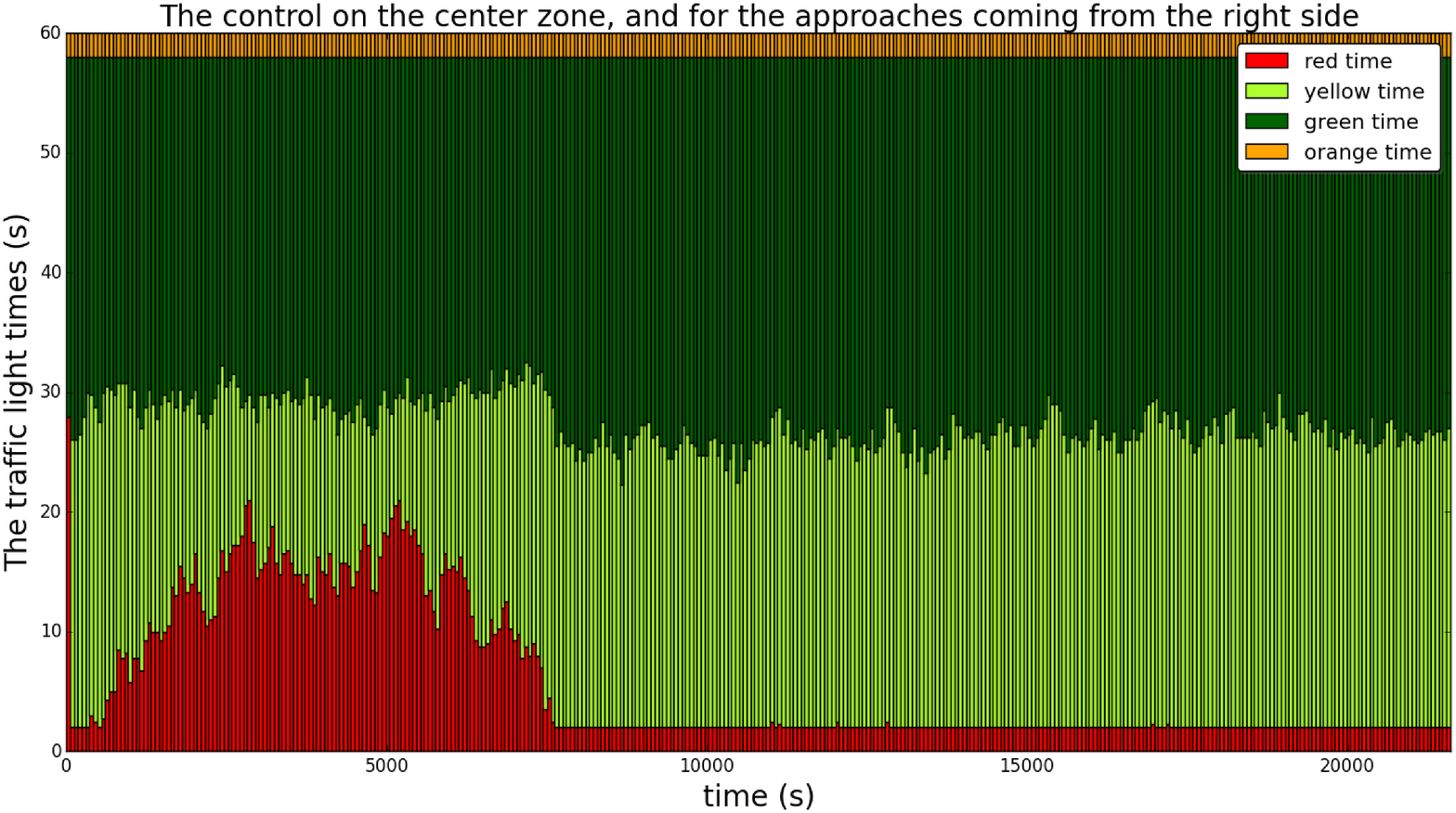} & \includegraphics[width=65mm]{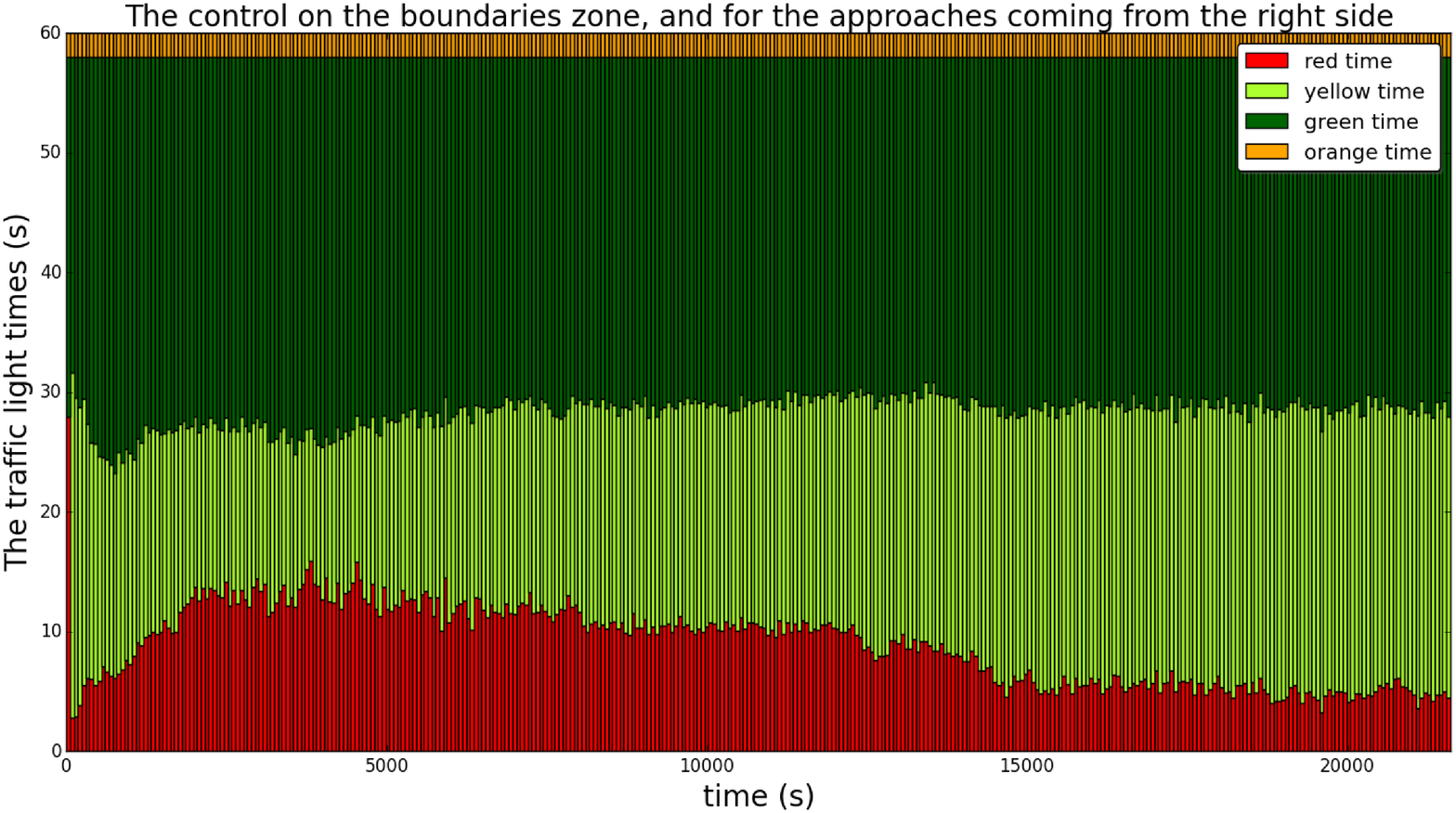}
\end{tabular}
\caption{Semi-decentralized TUC. The control in terms of the traffic light times into the cycle time, through the simulation time,
  on different zones (center and boundaries), and for approaches coming from left and right sides.}
 \label{fig-all}  
\end{center}
\end{figure}

In Figure \ref{fig-all}, we give the results of simulation for the semi-decentralized control.
We show on the first row the time-average number of vehicles in the circuits of the network.
On the second (resp. third) row of that figure, we show the control (in term of durations of the green, yellow and red times)
for the approaches coming from the left side (resp. right side) of the circuit junctions.
The left side column of the figure corresponds to the main circuit (the circuit of the central zone), while the right side column
corresponds to the secondary circuits (the circuits on the boundary of the network).

We observe on the first row of Figure \ref{fig-all} that the main circuit is more cleared out than the secondary circuits.
This observation confirms our intuition given above.
We see in the second and third rows of Figure \ref{fig-all} that the control frees the approaches coming from the left side
of the junctions' main circuit and limits the flow on the antagonistic approaches of the same circuit, while it does the opposite for the
secondary circuits.

Figure~\ref{fig-all} shows another important result, which is that the yellow time is almost fully used (i.e. the red time is almost zero)
in the case of free traffic flow, while the red time appears with important values in case of congestion.
This result is very important because it confirms the importance that the activation as well as the duration of the local control
(the contention time window with yellow times) are both controlled by the centralized control, which optimizes them in function of the
state of the traffic in the network.

\section{Preliminary conclusions}

We presented in this article a TUC-based approach for the control of urban traffic.
By defining a time contention window inside the time cycle, we introduced a little of decentralization of the control.
We have implemented and simulated the new control on a small American-like network. The traffic has been simulated
using the Simulation of Urban MObility tool while the control has been implemented with Python.
We are aware that we need more investigations in order to validate our assertions.
For that we will improve the implementation of our control by better managing the contention time window, in particular
using communication network simulators.
On this small network, we showed that our approach is effective in terms of many parameters including
the total network capacity as well as the average travel time.
Another important result we obtained is the confirmation that the centralized control optimizes the activation as well as
the duration of the decentralized control (the contention time window) in function of the state of the traffic in the network.

\bibliographystyle{elsarticle-harv}
\bibliography{mybib}

\end{document}